\documentclass[12pt]{amsart}
\usepackage[cp1251]{inputenc}
\usepackage[english,russian]{babel}
\usepackage{amsmath}
\usepackage{amssymb}
\usepackage{amsfonts}
\usepackage{graphicx}

\newtheorem{lemma}{Лемма}
\newtheorem{theorem}{Теорема}

\newtheorem{problem}{Задача}

\begin{document}

 УДК{517.984.50}

\title[ Оператор градиент дивергенции  в $\mathbf{L}_{2}(G)$ ]
{Оператор  градиент дивергенции  в подпространствах
$\mathbf{L}_{2}(G)$ }
\author{\bf {Р.С.~Сакс}}

\maketitle
 Аннотация. {\small Автор изучает структуру пространства
  $\mathbf{L}_{2}(G)$ вектор-функций, квадратично интегрируемых
      по ограниченной односвязной  области $G$ трехмерного    пространства 
     с гладкой границей 
    и роль операторов   градиента дивергенции
     и ротора в построении базисов в  подпространствах
     ${\mathcal{{A}}}$ и  ${\mathcal{{B}}}$.

     Доказана       само-сопряженность расширения $\mathcal{N}_d$
   оператора $\nabla\mathrm{div}$ в подпро-странство
   $\mathcal{A} _ {\gamma}\subset {\mathcal{{A}}}$ и базисность
   системы его   собственных функций. Выписаны явные
   формулы решения спектральной   задачи в шаре и условия разложимости
   вектор-функции   в ряд Фурье по собственным функциям
   градиента дивергенции. Изучена разрешимость краевой
   задачи:\newline
   $\nabla\mathrm{div}\,\mathbf{u}+\lambda\,\mathbf{u}=\mathbf{f}$\, в\, $G$, \,
       $ (\mathbf {n}\cdot\mathbf {u})|_{\Gamma}=g$ в пространствах
         Соболева $\mathbf{H}^{s}(G)$ порядка $s\geq 0$ и в подпространствах.

 Попутно изложены аналогичные результаты
         для оператора ротор и его симметричного расширения $S$ в $\mathcal{B}$.

         Эта работа
         есть продолжение исследований автора\,\cite{saUMJ13}--
         \cite{sa15b}  .}

\section {Введение и основные результаты}

\subsection {Основные пространства и операторы} В статье мы рассмат-риваем линейные
пространства над полем $\mathbb{C}$ комплексных чисел.

Через $\mathbf{L}_{2}(G)$ обозначаем пространство Лебега
вектор-функций, квадратично интегрируемых в $G$    с внутренним
произведением \newline $(\mathbf {u},\mathbf {v})= \int_G \,\mathbf
{u}\cdot\overline{\mathbf {v}}\,d\,\mathbf {x}$  и нормой
$\|\mathbf{u}\|~= (\mathbf {u},\mathbf {u})^{1/2}$.

  Пространство  Соболева порядка $s\geq 0$, состоящее из вектор-функций,
    принадлежащих $\mathbf{L}_{2}(G)$ вместе с обобщенными производными
  до порядка $s\geq 0$, обозначается через
$\mathbf{H}^{s}(G)$, $\|\mathbf {f}\|_s$ -норма его элемента
$\mathbf {f}$;   замыкание в $\mathbf{H}^{s}(G)$ пространства
$\mathcal{C}^{\infty}_0(G)$ обозначается через
$\mathbf{H}^{s}_0(G)$. Пространство Соболева отрицательного порядка
$\mathbf{H}^{-s}(G)$ двойственно к $\mathbf{H}^{s}_0(G)$
 (см. пространство $W_p^{(l)}(\Omega)$ при $p=2$ в $\S 3$ гл. 4 \cite{sob} ,
    $H^k(Q)$ в $\S 4$ гл. 3 \cite{mi}, а также   гл.1 в \cite{mag}).

 В области $G$  с гладкой границей  $\Gamma$ в каждой точке $y\in\Gamma$
  определена нормаль   $\mathbf {n}(y)$  к $\Gamma$.\,
Вектор-функция     $\mathbf {u}$ из $\mathbf{H}^{s+1}(G)$ имеет след
$ \gamma(\mathbf {n}\cdot\mathbf {u})$ на $\Gamma$ ее нормальной
компоненты, который принадлежит пространству Соболева-Слободецкого
$\mathbf{H}^{s+1/2}(G)$, $|\gamma(\mathbf {n}\cdot\mathbf
{u})|_{s+1/2}$- его норма.

{\it Свойства операторов ротор и  градиент дивергенции}

{\small Операторы градиент, ротор и   дивергенция определяются в
трехмерном векторном анализе \cite{zo}. Им соответствует оператор
$d$ внешнего диф-ференцирования на формах $\omega^k$ степени $k=0,1$
и 2. Соотношения
\newline
$dd\omega^k=0$ при $k=0,1$ имеют вид $\mathrm{rot}\,\nabla h=0$ и
$\mathrm{div}\,\mathrm{rot}\,\mathbf{u}=0$.

Формулы $\mathbf{u}\cdot\nabla
h+h\mathrm{div}\mathbf{u}=\mathrm{div}(h \,\mathbf{u}) $,\,
$\mathbf{u}\cdot\mathrm{rot}\,\mathbf{v}-
\mathrm{rot}\,\mathbf{u}\cdot\mathbf{v}=\mathrm{div}[\mathbf{v},\mathbf{u}]$,
где $[\mathbf{v},\mathbf{u}]$ - векторное произведение, и
интегрирование по области $G$ используются при определении
операторов $\mathrm{div}\,\mathbf{u}$  и $\mathbf{rot}\mathbf{u}$ в
$\mathbf{L}_{2}(G)$.}

 Пусть функция $h\in {H}^{1}(G)$, а
$\mathbf{u}=\nabla h$ - ее градиент. Через ${\mathcal{{A}}}(G)$
обозначим пространство ${\mathcal {{A}} } (G)=\{\nabla h, h\in
H^1(G)\}$, оно содержится в $\mathbf{L}_{2}(G)$ и содержит
подпространство
\begin{equation}\label{ag  1}\mathcal{A}_{\gamma} =
     \{\nabla h, h\in H^1(G),  \gamma(\mathbf {n}\cdot {\nabla\,h})=0
     \},\end{equation}
 $\mathcal{A} _ {\gamma}$ плотно в ${\mathcal {{A}} }(G)$ и содержит
 подпространство
 \begin{equation}\label{ag  1}\mathcal{A}^{2}_ {\gamma}=\{\mathbf {u}\in\mathcal{A} _ {\gamma}(G):
 \nabla div\mathbf {u}\in\mathcal{A} _ {\gamma}(G)\}.\end{equation}
 Оказывается, что  оператор градиент дивергенции
  имеет в $\mathcal{A}_{\gamma}$   само-сопряженное 
  расширение, а
именно, оператор 
  $\mathcal{N}_d:\mathcal{A}_{\gamma}(G)
 \longrightarrow \mathcal{A}_{\gamma}(G)$  с областью определения
 $\mathcal{A}^{2}_ {\gamma}$,
  совпадающий с $\nabla div$ на  $\mathcal{A}^{2}_ {\gamma}
  \subset \mathbf{H}^2(G)$,
  --самосопряжен.   Обратный оператор
 $\mathcal{N}^{-1}_d:\mathcal{A} _ {\gamma}\rightarrow \mathcal{A}^{2}_ {\gamma}$
  вполне непрерывен (теорема 2).

   Следовательно,   этот оператор имеет
полную систему собственных функций, отвечающих ненулевым собственным
значениям:
\[-\nabla\mathrm{div}\mathbf{q}_{j}=\mu_j \mathbf{q}_{j},\quad \mu_j\in M
\subset \mathbb{R},\]
\[\mathbf{a}(x)=\sum_{\mu_j\in M}(\mathbf{a},\mathbf{q}_{j})\mathbf{q}_{j}
(x),\quad \text{если}\quad \mathbf{a}(x)\in{\mathcal {{A}} _
{\gamma}} (G),\quad \|\mathbf{q}_{j}\|=1,\]
  Ортогональное дополнение ${\mathcal {{A}} } (G)$ в
$\mathbf{L}_{2}(G)$ обозначим через $\mathcal{B}(G)$.

 Пространство ${\mathcal{{B}}}(G)$  в
   \cite{giyo}  Z.Yoshida и Y.Giga обозначают как
  ${L}_{\sigma}^2(G)$, А.Фурсиков в \cite{afu}-- как $\mathbf {V}(G)$.
   В обобщенном смысле оно формулируется так:
  \[\mathcal{B}(G)=\{\mathbf{u}\in\mathbf{L}_{2}(G):
   \,\mathrm{div}\,\mathbf{u}=0,\,\mathbf{n}\cdot \mathbf{u}|_{\Gamma}=0 \}.\]
  Если граница
области $G$ имеет положительный род $\rho$, то $\mathcal{B}(G)$
содержит  конечномерное подпространство
\begin{equation}\label{bh  1} \mathcal{B}_H=\{\mathbf{u}\in\mathbf{L}_{2}(G):
\mathrm{rot}\,\mathbf{u}=0,\,\,\mathrm{div}\,\mathbf{u}=0,
\,\,\mathbf{n}\cdot \mathbf{u}|_{\Gamma}=0\}. \end{equation} Его
размерность
 равна $\rho$ \cite{boso}, а базисные функции
 $h_j\in \mathcal{C}^\infty(G)$. Гладкость обобщенных рещений системы
  \eqref{bh  1}
 доказал  Г.Вейль \cite{hw}, 1941, а решений полигармонических уравнений
 С.Соболев \cite{so37}, 1937.

   Ортогональное дополнение $\mathcal{B}_H$ в
$\mathcal{B}(G)$ обозначим через $\mathbf {V} ^{0} (G)$. Значит,
\begin{equation}\label{bhr  1}
\mathbf{L}_{2}(G)= \mathcal{A}(G)\oplus\mathcal{B}(G),\quad
{\mathcal{B}(G)= \mathbf {V} ^{0} (G)\oplus\mathcal{B}_H}
(G).\footnote {{\small 
См.  Z.Yoshida и
Y.Giga  \cite{giyo}. Разложение 
$\mathbf{L}_{2}(\Omega)$ приводят также Г.Вейль \, \cite{hw},
С.Л.Соболев \cite{so54}, О.А.Ладыженская 
   \cite{lad}, К.Фридрихс \cite{fri},  Н.Е.Кочин,  И.А.Кибель, Н.В.Розе
     \cite{kkr} и
   Э.Быховский и Н.Смирнов \cite{bs}.  }}
\end{equation} $\mathbf {V} ^{0} (G)$ содержит подпространство
\[\mathbf{W}^{1}=\{\mathbf {u}\in\mathbf {V}^0(G):\mathrm{rot}\mathbf
{u}\in\mathbf {V}^0(G)\}.\] Z.Yoshida и Y.Giga \cite{giyo} показали,
что оператор ротор имеет в $\mathbf{V}^{0}(G)$ самосопряженное
расширение: оператор  $S:\mathbf{V}^{0}\rightarrow \mathbf{V}^{0}$ с
областью определения  $\mathbf{W}^{1}$, совпадающий с
  $\mathrm{rot}\, \mathbf{u}$ на 
  $\mathbf{W}^{1}\subset\mathbf {H}^1(G)$ -- самосопряжен.

Собственные функции $\mathbf{u}_j$ ротора, отвечающие ненулевым
собственным значениям, образуют базис в $\mathbf{V}^{0}(G)$.

 В гидродинамической интерпретации \cite{kkr} им соответствуют потоки, имеющие
 ненулевую завихренность, а собственным функциям $\mathbf{q}_j$ оператора
 градиент дивергенции--потенциальные (безвихревые) потоки.

 Приложения смотрите в работах \cite{koz, tay, motuva, cdtgt}, а также
 \cite{saTMF, saha, saUMJ, saUMJ13}.

 {\small  В  шаре $B$ радиуса $R$ собственные функции
 $\mathbf{u}^{\pm}_{\kappa}$ ротора, отвечающие
  ненулевым собственным значениям
  $\pm\lambda_{\kappa}=\pm\rho_{n,m}/R$ и собственные функции
  $\mathbf{q}_{\kappa}$
градиента дивергенции с собственными значениями $-\nu_{\kappa}^2$,
$\nu_{\kappa}=\alpha_{n,m}/R,$ выражаются  явными формулами
\cite{saUMJ13}, причем
    \[\mathbf{rot}\,\mathbf{u}^{\pm}_{\kappa}=\pm\lambda_{\kappa}\,
\mathbf{u}^{\pm}_{\kappa}, \quad
\nabla\,div\,\mathbf{u}^{\pm}_{\kappa}=0,\quad
\gamma\mathbf{n}\cdot\mathbf{u}^{\pm}_{\kappa}=0;
 \quad \kappa=(n,m,k),\,\]
\[ \mathbf{rot}\,\mathbf{q}_{\kappa}=0, \quad
-\nabla\,div\,\mathbf{q}_{\kappa}=\nu_{\kappa}^2\mathbf{q}_{\kappa},
\quad \gamma\mathbf{n}\cdot\mathbf{q}_{\kappa}=0, \quad
 \quad\,\, |k|\leq n ,\] где
числа $\pm\rho_{n,m}$ и $\alpha_{n,m}$ - нули функций $\psi_n$ и их
производных $\psi_n'$, а
 \begin{equation}
\label{psi__1_}\psi_n(z)=(-z)^n\left(\frac{d}{zdz}\right)^n\frac{\sin
z}z, \quad n\geq 0,\,\,m\in \mathbb{N}.\end{equation}

Cобственные функции каждого из операторов взаимно ортогональны и их
совокупная система полна в ${\mathbf{{L}}_{2}}(B)$ \cite{saUMJ13}.

Найдены необходимые и достаточные условия на вектор-функцию
$\mathbf{v}$ из $\mathcal{A}_{\gamma}(B)$, при которых ее ряд Фурье
по собственным функциям оператора $\nabla\,div$ сходится в норме
пространства  $\mathbf{H}^{s}(B)$ порядка $s>0$ . Эти условия
состоят в  принадлежности $\mathbf{v}$ пространству
$\mathcal{A}_{\mathcal{K}}^{s}(B)$ (п.3.6 Теорема 3).}

 \subsection{ Структура работы }
 В $\S 2$
 мы изучаем в
ограниченной области $G$ с гладкой границей $\Gamma$ разрешимость
краевой задачи  $\gamma\mathbf{n}\cdot\mathbf{u}=g$
 для системы
$\nabla\mathbf{div}\mathbf{u}+\lambda\mathbf{u}=\mathbf{f}$ в  $G$ в
пространствах $\mathbf{H}^{s}(G)$ при  $s\geq 2$.

Эта система не является эллиптической \cite{sa75}, но при
$\lambda\neq 0$ принадлежит классу Б.Вайберга и В.Грушина
\cite{vagr} систем, приводимых  к  эллиптическим: ее расширение
является эллиптической 
системой.

 Мы доказываем, что указанная краевая задача удовлетворяют условиям
эллиптичности из Теоремы 1.1 Солонникова в \cite{so71}.

Следовательно, оператор $\mathbb{B}$ задачи в пространствах (\ref{op
2}) имеет левый регуляризатор, конечномерное ядро  и выполняются
априорная оценка:
\begin{equation} \label{apgro s} C_s\|\mathbf{u}\|_{s+2}
\leq|\lambda|\,\|\mathrm{rot}^2\,\mathbf{u}\|_{s}+
\|\nabla\mathrm{div}\,\mathbf{u}\|_{s}+
|\gamma({\mathbf{n}}\cdot\mathbf{u})|_{s+3/2}+ \|\mathbf{u}\|_{s}
 \end{equation}
 Разрешимость этой задачи зависит от пространств, к
которым принадлежат $\mathbf{f}$ и $g$ ( пп.2.3  ). 
 Пространство $\mathcal{B}(G)$ 
 принадлежит ядру оператора градиента дивергенции в
$\mathbf{L}_{2}(G)$, поэтому уравнение
$\nabla\mathbf{div}\mathbf{u}+\lambda\mathbf{u}=\mathbf{g}$ при
$\mathbf{g}\in \mathcal{B}(G)$ сводится к уравнению
$\lambda\mathbf{u}=\mathbf{g}$.

  На подпространство $\mathcal{A}_{\gamma}$
  он продолжается как самосопряженный оператор
$\mathcal{N}_d + \lambda\,I: \mathcal{A}_{\gamma}\rightarrow
\mathcal{A}_{\gamma} $. Найдены необходимые и достаточные условия
обратимости этого оператора (Теорема 2).

В  $\S\, 3$ спектральная задача для оператора градиент дивергенции в
  области с гладкой границей сводится к решению спектральной
задачи Неймана для скалярного оператора Лапласа.

 В шаре ее решения вычислены  явно\,\cite{vla}. В результате мы получаем формулы
\eqref{qom 1} собственных функций $\mathbf{q}_{\mu}(\mathbf{x})$
 градиента дивергенции.

Они являются также собственными функциями ротора  с {\it нулевым}
собственным значением.

   \section{ Градиент дивергенции
 в ограниченной области}

\subsection{Краевая задача:} Пусть  $G$ - ограниченная область
  в ${{R}^{3}}$ с гладкой
границей $\Gamma $, $\mathbf{n}$- внешняя нормаль к $\Gamma $.

В частности, $G$ может быть шаром $B$,  $|x|<R$,
 с границей $S$.
\begin{problem}В  области $G$ и на ее границе $\Gamma$
 заданы векторная и скалярная функции $\mathbf{f}$ и ${g}$, найти
вектор-функцию $\mathbf{u}$, такую что
\begin{equation}
\label{grd__2_}
 \nabla\, \text{div}\mathbf{u}+\lambda\, \mathbf{u}=\mathbf{f}
 \quad  \text{в}\quad G,
\end{equation}
\begin{equation}
\label{Gr__1_}\mathbf{n}\cdot \mathbf{u}{{|}_{\Gamma }}=g,
\end{equation}
 где $\mathbf{n}\cdot \mathbf{u}$ - скалярное произведение векторов
$\mathbf{u}$ и $\mathbf{n}$.\end{problem}

Эта {\it задача не эллиптична}. Оператор $\nabla\,
\text{div}+\lambda\text{I}$ второго порядка не является
эллиптическим, так как ранг его символической матрицы $\nabla\,
\text{div}(i\xi)$  равен единице при всех $\xi\in
\mathcal{R}^3\backslash 0$ и меньше трех \cite{sa75}.

Б.Вайнберг и В.Грушин   \cite{vagr} 1967 определили на гладком
многообразии $X$ без края  класс {\it равномерно неэллиптических
систем } (РНС) cингулярных интегро-дифференциальных уравнений и
класс матричных с.и.д.операторов, {\it глобально приводимых к
эллиптическим матрицам}, и доказали их эквивалентность. Эти
определения требуют введения дополнительных понятий.

 Мы приведем их  для
систем дифференциальных уравнений с постоянными коэффициентами,
который обозначим как (РНСp)

 {\small Система дифференциальных уравнений, $S(D)u=f$ порядка $m$, из этого
класса обладает свойствами: \newline
 а){\it ее символическая матрица
$S_0(i\xi)$ имеет постоянный ранг при всех}
$\xi\in\mathcal{R}^3\backslash 0$.  Это позволяет построить
аннулятор $C(D)$ оператора $S_0(D)$ такой, что
$(CS_0)(D)\mathbf{u}\equiv 0$ на $X$ и определить
\newline б){\it расширенную систему \quad $Su=f,// CSu=Cf$
порядка~$m//l$.}
\newline
Ее символическая матрица $S_0(i\xi),//(CS)_0(i\xi)$ определяется
младшей частью оператора $S(D)$ и дополняет матрицу $S_0(i\xi)$.
\newline в){\it Если ранг  расширенной матрицы максимален, то
исходная система $Su=f$ принадлежит классу (РНС1) и  степень ее
приводимости равна единице.}

г){\it Если система $Su=f$ такова, что ранг  расширенной матрицы не
максимален, но постоянный, то процесс повторяется} и при
определенных условиях система принадлежит классу (РНС2). И  так
далее.

Б.Вайнберг и В.Грушин  \cite{vagr} доказали, что на замкнутом
многооббразии $X$ система $Su=f$ класса (РНСp) являются разрешимой
по Фредгольму или Нетеру в пространствах Соболева $\mathbf{H}^s(X)$,
если $f\in \mathbf{H}^{s-m+p}(X)$, где $s\geq m$ целое. В качестве
примера оператора из класса (РНС1)  приводится оператор $d+\ast$ на
дифференциальных формах степени $k$ на $2k+1$-мерном многообразии
$X$.}

  Покажем, что  система \eqref{grd__2_} 
  принадлежит классу (RNS1) при  $\lambda\ne 0$. Действительно,
  оператор $\nabla\, \text{div}+\lambda\text{I}$ таков, что \newline
   а)  {\it ранг его символической матрицы
$\nabla\, \text{div}(i\xi)$ равен
единице  при любых $\xi\neq0$}, 
\newline
 б)оператор $\nabla\,\text{div}$  имеет  левый
 аннулятор $\text{rot}$, так как $ \text{rot}\,\nabla\, \text{div} \mathbf{u}=0$
 для любой $\mathbf{u}\in \mathcal{C}^3(G)$,\newline
 в) {\it $(6\times3)$-оператор
 $\nabla\,\text{div}$// $\lambda\text{rot}$ порядка $2//1$
 эллиптичен.}

 Его символическая матрица имеет максимальный ранг (=3)
  при выборе определенных 
  порядков $s_k$ и $t_j$ для его строк и столбцов, а именно при
$s_k=0$ для $k=1,-,3$ и $s_k=-1$ для $k=4,-,6$; и при $t_j=2$ для
$j=1,-,3$\, \cite{sa75}. Ранг матрицы $\nabla\,\text{div}(i\xi )$//
$\lambda\text{rot}(i\xi $ равен  трем при всех $\xi \neq 0$,
 поэтому расширенная система:
  \begin{equation}
\label{rodi__2_} \nabla\, \text{div}\,\mathbf{u}+\lambda
\mathbf{u}=\mathbf{f},\quad \lambda \text{rot}\,\mathbf{
u}=\text{rot}\,\mathbf{f}\end{equation} является эллиптической по
Даглису-Ниренбергу. В общем случае  (при  $\mathbf{F}$ вместо 
$\text{rot}\,\mathbf{f}$)эта система переопределена.

 Оказывается, что формально переопределенная краевая задача \eqref{Gr__1_}, \eqref{rodi__2_},
  эллиптична по определению В.А.Солонникова \cite{so71}. Это означает что

1)  система \eqref{rodi__2_} эллиптична,

2) граничный оператор $\gamma\mathbf{n}\cdot\mathbf{u}$ "накрывает"
\,оператор системы \eqref{rodi__2_}.

Первое утверждение выполняется, если 

 $1^0)$ однородная система линейных
алгебраических уравнений:
\begin{equation}
\label{cdx  2} \lambda\,\text{rot}(i\xi )\mathbf{w}=0, \quad
(\nabla\,\text{div})(i\xi )\mathbf{w}=0, \quad \forall \xi\neq 0
\end{equation}
 c параметром $\xi \neq 0$ имеет только тривиальное
 решение $ \mathbf{w}(\xi)=0$;

Пусть  $\tau$ и $\mathbf{n}$- касательный и нормальный векторы к
$\Gamma$ в точке $y\in \Gamma$ и $|\mathbf{n}|=1$.
  Второе утверждение выполняется, если

  $2^0)$   однородная система линейных
дифференциальных уравнений ( на полуоси\, $z\geq 0$ с параметром
$|\tau| > 0$):
\begin{equation}\label{cdz  3}
\lambda\text{rot}(i\tau+\mathbf{n} d/dz ) \mathbf{v}=0,\quad
(\nabla\text{div})(i\tau+\mathbf{n} d/dz )\mathbf{v}=0 \,\, 
\end{equation}  с условиями:
$\mathbf{n}\cdot \mathbf{v}|_{ z=0}=0$ и $\mathbf{v}\rightarrow 0$
при $z\rightarrow + \infty$ имеет только тривиальное решение
$\mathbf{v}(y,\tau; z )$.

Для доказательства $1^0)$, $2^0)$ воспользуемся соотношением
\begin{equation}\label{cc  2} \text{rot}\, \text{rot}\,\mathbf{v}=
 -\Delta \mathbf{v}
 + \nabla \text{div} \mathbf{v}.\end{equation}
 $1^0).$  Из уравнений  (\ref{cdx  2}) вытекает уравнение
$-\Delta(i\xi)\mathbf{w}=0$, которое распадается на три скалярных
уравнения $|\xi | ^2 w_j(\xi) =0$. Значит, $\mathbf{w}=0$ ибо $|\xi
|\neq 0$, и
 система \eqref{rodi__2_} -- эллиптична.
\newline
 $2^0).$  Из уравнений (\ref{cdz  3}) вытекает уравнение $ (-|\tau| ^2 +
(d/dz)^2)\mathbf{v} = 0$ с параметром $|\tau |> 0$. Следовательно,
$\mathbf{v}=\mathbf{w} e^{ -|\tau|z}$. Эта вектор-функция
удовлетворяет уравнениям (\ref{cdz 3}), если вектор
 $\mathbf{w}$ таков, что
   $\quad \omega\times \mathbf{w}=0\quad \omega ({\omega}' \cdot
\mathbf{w})=0$, где $\omega= i\tau -| \tau| \mathbf{n}$ -
вектор-столбец. Так как  $\overline{{\omega}'} \cdot
\omega=|\tau|^2\neq 0$, то ${\omega}' \cdot \mathbf{w}=0$.

Уравнения $\omega \times \omega=0,$ \,
 ${\omega}'\cdot\omega=0$, имеют решение
 $\mathbf{w}=c\,\omega$, где $c$ - постоянная.
  Граничное условие приводит нас к равенству
  $|\tau |c=0$.  Следовательно,   $c=0$ ибо $|\tau| > 0$ и $\mathbf{v}=0$.

Итак, задача \eqref{Gr__1_},  \eqref{rodi__2_}  эллиптична по
Солонникову .

 Мы скажем в этом случае, что краевая задача
 \eqref{grd__2_}, \eqref{Gr__1_} при $\lambda\neq 0$
  является обобщенно эллиптической.

\subsection{Оператор задачи 1 в пространствах Соболева}
Пусть $\mathbf{u}$ принадлежит пространству ${\mathbf{H}^{s+2}}(G)$,
то-есть каждая компонента $u_j\in H^{s+2}(G)$. Тогда
$\nabla\text{div} \mathbf{u}$ принадлежит $\mathbf{H}^{s}(G)$, и
$\text{rot}^{2}\mathbf{u}$ принадлежит $\mathbf{H}^{s}(G)$. Поэтому
вектор-функция $\mathbf{f}:=\nabla\text{div} \mathbf{u}+ \lambda
\mathbf{u}$ принадлежит пространству
\begin{equation}\label{pr  2} {\bf{F}^{s}}(G)=\{\mathbf{f}\in
{\mathbf{H}^{s}(G)}: \text{rot}^2\,\mathbf{f}\in
 \mathbf{H}^{s}(G)\},\end{equation} которое снабдим  нормой
 \[\|\mathbf{v}\|_{\mathbf{F}^{s}}=(
 \|\mathbf{v}\|^2_{s}+
 \|\text{rot}^2\mathbf{v}\|^2_{s})^{1/2}.\] Функция
$g:=\gamma({\mathbf{n}}\cdot\mathbf{u})\equiv\mathbf{n}
\cdot\mathbf{u}|_{\Gamma}$ принадлежит пространству
$H^{s+3/2}(\Gamma)$. Следовательно, при $\lambda\neq 0$ задаче
соответствует ограниченный оператор
\begin{equation}\label{op  2} \mathbb{B}\mathbf{u}\equiv \begin{matrix}
\nabla\text{div}\,\mathbf{u}+\lambda\,\mathbf{u} \\
\gamma\mathbf{n}\cdot \mathbf{u}\end{matrix}:
\bf{H}^{s+2}(G)\rightarrow
\begin{matrix}\bf{F}^{s}(G)\\ H^{s+3/2}(\Gamma)\end{matrix}.
\end{equation}
Согласно Теореме 1.1 из работы Солонникова \cite{so71}, о
переопределенных эллиптических краевых задачах,  в ограниченной
области $G$   с гладкой границей $\Gamma\in \mathcal{C}^{s+2} $,
обобщенно эллиптический оператор \eqref{op 2} имеет левый
регуляризатор, то-есть
 ограниченный оператор  $\mathbb{B}^L$ такой, что
$\mathbb{B}^L\mathbb{B}=\mathbb{I}+\mathbb{T}$, где $\mathbb{I}$ -
единичный, а $\mathbb{T}$ - вполне непрерывный операторы, область
значений замкнута, и существует постоянная $C_s >0$ такая, что
выполняется оценка:
\begin{equation} \label{apgro s} C_s\|\mathbf{u}\|_{s+2}
\leq|\lambda|\,\|\mathrm{rot}^2\,\mathbf{u}\|_{s}+
\|\nabla\mathrm{div}\,\mathbf{u}\|_{s}+
|\gamma({\mathbf{n}}\cdot\mathbf{u})|_{s+3/2}+ \|\mathbf{u}\|_{s}
 \end{equation}
где $\|\mathbf{u}\|_{s+2}$ норма  $\mathbf{u}$ в
$\mathbf{H}^{s+2}(G)$,
 $|\gamma(\mathbf{n}\cdot\mathbf{u})|_{s+3/2}$ --  норма следа
нормальной компоненты $\mathbf{u}$ на $\Gamma $
 в
${H}^{s+3/2}(\Gamma)$, $s\geq 0$\,  \cite{so71,sa75}. \,
 Итак, имеет место
\begin{theorem}
Оператор $\mathbb{B}$ в пространствах (\ref{op  2}) имеет левый
регуляризатор. Его ядро $\mathcal{M}$ конечномерно, область значений
замкнута и выполняется априорная оценка \eqref{apgro s}.
 \end{theorem}

Из этой теоремы и  оценки следует, что при $\lambda \ne 0$

 a){\it число линейно независимых решений однородной
задачи 1 конечно,}

 b){\it любое ее обобщенное решение 
 бесконечно дифференцируемо вплоть до
границы, если граница области бесконечно дифференцируема.}

Замечание. Оценку \eqref{apgro s} я не видел у других авторов.
Известна \cite{giyo} оценка:  существует постоянная $C_s >0$ такая,
что 
\begin{equation} \label{apro s} C_s\|\mathbf{u}\|_{s+1}
\leq\|\mathrm{rot}\,\mathbf{u}\|_{s}+
\|\mathrm{div}\,\mathbf{u}\|_{s}+
|\gamma({\mathbf{n}}\cdot\mathbf{u})|_{s+1/2}+ \|\mathbf{u}\|_{s}
 \end{equation}

\subsection{Оператор $\nabla\mathrm{div}+\lambda I$ в подпространствах }
На подпространстве $\mathcal{B}$ в $\mathbf{L}_2(G)$, ортогональном
  $\mathcal{A}=\{\mathbf{u}=\nabla\, h:\,\, h\in
H^1(G)\}$
 оператор $\nabla\, \text{div}\mathbf{u} +\lambda \mathbf{u} $
является алгебраическим оператором вида: $\lambda \mathbf{u} $.

 Пространство
$\mathcal{A}_{\gamma}=\{\mathbf{u}=\nabla\, h:\,\, h\in H^1(G),\quad
(\mathbf{n}\cdot\mathbf{u})|_{\Gamma}=0 \}$
плотно в $\mathcal{A}$, так как функции из $\mathcal{C}_0^\infty \cap
\mathcal{A}_{\gamma}$ плотны в $\mathbf{L}_2(G)$. Пространство
  \[\mathcal{A}^2_ {\gamma}(G)=\{ \mathbf{v}\in
\mathcal{A}_{\gamma}: \nabla\mathrm{div} \mathbf{v}\in
\mathcal{A}_{\gamma}\}.\]
 плотно в
$\mathcal{A}_{\gamma}$ и
 содержится в $\mathbf{H}^2(G)$  согласно оценке \eqref{apgro s}.

Введем оператор $\mathcal{N}_d  :\mathcal{A}_{\gamma}
 \rightarrow \mathcal{A}_{\gamma}$ с  областью
определения   $\mathcal{A}^2_ {\gamma}(G)$,  который при
$\mathbf{u}\in\mathcal{A}^2_ {\gamma}(G)$ 
 совпадает с $\nabla\mathrm{div}\mathbf{u}$.

Оператор $\mathcal{N}_d+\lambda I$, где $I$ - единичный оператор,
  является самосопряженным (эрмитовым\, \cite{vla}).
 Действительно, по формуле Гаусса-Остроградского
\begin{equation}
\label{gri__2_} \int_G(\nabla\mathrm{div}\mathbf{u}+\lambda
\mathbf{u})\cdot\,\mathbf{v} dx=
\int_G\mathbf{u}\cdot\,(\nabla\mathrm{div}\mathbf{v}+\lambda
\mathbf{v})dx+\end{equation}
\[\int_{\Gamma}[(\mathbf{n}\cdot\mathbf{v})\mathrm{div}\mathbf{u}+
(\mathbf{n}\cdot\mathbf{u})\mathrm{div}\mathbf{v}] \, d S.\]

Если вектор-функции $\mathbf{u}$ и $\mathbf{v}$ принадлежат
$\mathcal{A}^2_ {\gamma}(G)$, то граничные интегралы пропадают,
остальные интегралы сходятся. Следовательно,
\begin{equation}
\label{gri__2_} ((\nabla\mathrm{div}\mathbf{u}+\lambda
\mathbf{u}),\mathbf{v})=
(\mathbf{u},(\nabla\mathrm{div}\mathbf{v}+\lambda \mathbf{v})) \quad
\text{в}\quad \mathbf{L}_2(G).\end{equation}

\begin{lemma} Область значений $\mathcal{R}(\mathcal{N}_d)$
оператора $\mathcal{N}_d$ совпадает с
$\mathcal{A}_{\gamma}$.\end{lemma} Доказательство. Из определения
$\mathcal{D}(\mathcal{N}_d)$ видим, что
$\mathcal{R}(\mathcal{N}_d)\subset \mathcal{A}_\gamma$. Положим
\[ \widehat{A}=\{\mathbf{v}=P_a\,\mathbf{w};\, \mathbf{w}\in H^2(G):
(\mathbf{n}\cdot\nabla)\,div\, \mathbf{w})|_\Gamma=0,\}\] где $P_a$
есть ортогональный проектор из $\mathbf{L}_2(G)$ на
$\mathcal{A}_\gamma$. 
 По определению $\mathbf{w}-P_a\,\mathbf{w}
\,\in Ker(\nabla\,div)$ , значит
$\widehat{A}\subset\,\mathcal{D}(\mathcal{N}_d)$, потому что для
$\mathbf{v}=P_a\,\mathbf{w}\in \widehat{A}$, имеем
$\mathbf{v}\in\mathcal{A}_\gamma$ и
$(\nabla\,div)\,\mathbf{v}=(\nabla\,div)\,\mathbf{w}\in
\mathcal{A}_\gamma$.

Далее,\, $\mathcal{N}_d(\widehat{A})$=$\{(\mathcal{N}_d)\mathbf{v},
\,\mathbf{v}\in
\widehat{A}\}$=
$$\{\nabla\,div\,P_a\,\mathbf{w}=\nabla\,div\,\mathbf{w},
\,\, \mathbf{w}\in H^2(G):(\mathbf{n}\cdot
\nabla)\,div\,\mathbf{w}|_\Gamma=0\}.$$ Значит,
$\mathcal{A}_\gamma\subset\mathcal{R}(\mathcal{N}_d)$.
Следовательно,\, $\mathcal{A}_\gamma=\mathcal{R}(\mathcal{N}_d)$.
Лемма доказана.

Так как $\mathcal{A}_\gamma$ ортогонально $Ker\mathcal{N}_d$,
оператор $\mathcal{N}_d$ имеет единственный обратный
$\mathcal{N}_d^{-1}$ определенный на $\mathcal{A}_\gamma$. Оператор
 $\mathcal{N}_d^{-1}:\mathcal{A}_{\gamma}\rightarrow
 \mathcal{A}^2_\gamma(G)$ имеет точечный спектр, который
 не содержит точек накопления кроме нуля.

В следующем параграфе мы покажем, что
  собственные значения оператора $-\mathcal{N}_d$ совпадают
  с собственными значениями оператора Лапласа-Неймана
  $\mathcal{N}_{(-\Delta)}$ \eqref{gen   1} и являются положительными.
 Их множество счетно и каждое из собственных значений $\mu$ имеет
  конечную кратность. Перенумеруем их в порядке возрастания:
  $0<\mu_1\leq \mu_2\leq ...$, повторяя $ \mu_k$ столько раз, какова
  его кратность. Соответствующие собственные функции обозначим через
$\mathbf{q}_{1}, \mathbf{q}_{2}$, ...,  так чтобы каждому
собственному значению $\mu_{k}$ соответствовала только одна
собственная функция $\mathbf{q}_{k}$:
\[-\mathcal{N}_d\,\mathbf{q}_{k}=\mu_k\,\mathbf{q}_{k},\quad
k=1,2,...\] Собственные функции, соответствующие одному и тому же
собственному значению, выберем ортонормальными, используя процесс
ортогонализации Шмидта \cite{vla}. Собственные функции,
соответствующие различным собственным значениям, ортогональны. Их
нормируем.

Таким образом, в пространстве $\mathcal{A}_{\gamma}$ оператор
$-\mathcal{N}_d$ позволяет построить ортонормальный базис
$\{\mathbf{q}_{j}(\mathbf{x})\}$, такой что
\begin{equation}
\label{spr__4_}\mathbf{f}(\mathbf{x})=\sum_{j=1}^\infty
(\mathbf{f},\mathbf{q}_{j})\mathbf{q}_{j}(\mathbf{x}) \end{equation}
и ряд сходится в $\mathbf{L}_{2}(G)$ для любой  $\mathbf{f}\in
\mathcal{A}_{\gamma}(G)$.

 Так как
 $(\mathcal{N}_d+\lambda\,I)\mathbf{f}\in \mathcal{A}_{\gamma}(G)$
 по определению, то

\begin{equation}
\label{sp__2_} (\mathcal{N}_d+\lambda\,I)\mathbf{f}=
\sum_{j=1}^\infty[(\lambda-\mu_j)(\mathbf{f},\mathbf{q}_{j})
\mathbf{q}_{j}]
\end{equation} и ряд сходится в $\mathbf{L}_{2}(G)$.
Если $\lambda=\mu_{j_0}$,
 то соответствующее слагаемое в этом ряду исчезает.

Если элемент $(\mathcal{N}_d+\lambda\,I)^{-1}\mathbf{f}\in
\mathcal{A}_{\gamma}(G)$, то
\begin{equation}
\label{sp__4_}
(\mathcal{N}_d+\lambda\,I)^{-1}\mathbf{f}=
\sum_{j=1}^\infty[(\lambda-\mu_j)^{-1}(\mathbf{f},\mathbf{q}_{j})\mathbf{q}_{j}
]
\end{equation} и ни одно из слагаемых этого ряда не обращается в
бесконечность. Это означает, что $(\mathbf{f},\mathbf{q}_{j})=0$ при
$\lambda=\mu_j=\mu_{j_0}$, то-есть функция $\mathbf{f}$ ортогональна
всем собственным функциям $\mathbf{q}_{j}(\mathbf{x})$ градиента
дивергенции, отвечающим собственному значению $\mu_{j_0}$.
 Итак, имеет место
\begin{theorem}
a). Оператор $\mathcal{N}_d:\mathcal{A}_{\gamma}\rightarrow
\mathcal{A}_{\gamma}$   является
 самосопряженным.
  Его спектр  $\sigma(\mathcal{N}_d )$
 точечный  и действительный.
  Семейство собственных функций $\mathbf{q}_{j}(x)$ оператора
  $\mathcal{N}_d$ образует полный ортонормированный базис в  пространстве
   $\mathcal{A}_{\gamma}$;
    разложение   $\mathbf{a}(x)\in{\mathcal {{A}} _{\gamma}} (G)$ имеет вид
\begin{equation}
\label{spr__5_}\mathbf{a}(x)=\sum_{j=1}^\infty(\mathbf{a},\mathbf{q}_{j})
\mathbf{q}_{j}(x),\quad
\|\mathbf{q}_{j}\|=1.\end{equation}
  b). Если $\lambda$ не совпадает ни с одним из
   собственных   значений оператора $\mathcal{N}_d$, то
   оператор $\mathcal{N}_d+\lambda\,I:
    \mathcal{A}_{\gamma}\rightarrow\mathcal{A}_{\gamma}$ однозначно обратим,
    и его обратный   задается формулой   \eqref{sp__4_}.
Если $\lambda=\mu_{j_0}$, то он обратим тогда и только
  тогда,когда
  \begin{equation}
\label{urz _1_}\int_G \mathbf{f}\cdot \mathbf{q_j}\, dx=0\quad
\text{для}\,\,\forall \mathbf{q_j}: 
\mu_j=\mu_{j_0}.
\end{equation}
Ядро оператора $\mathcal{N}_d+\mu_{j_0}\,I$ определяется
собственными функциями $\mathbf{q_j}(\mathbf{x})$, собственные
значения которых равны $\mu_{j_0}$:
\begin{equation} \label{ker__1_}
Ker(\mathcal{N}_d+\mu_{j_0}\,I)= \sum_{\mu_j=\mu_{j_0}}
c_j\,\mathbf{q}_{j}(\mathbf{x}), \quad \text{для}\,\,\forall c_j \in
\mathcal{R}.\end{equation}
  \end{theorem}
Согласно Лемме 2 $\S 3$ собственные значения оператора
$-\mathcal{N}_d$ совпадают  с собственными значениями оператора
Лапласа-Неймана  $\mathcal{N}_{(-\Delta)}$ \eqref{gen   1} и
являются положительными. Значит, оператор $\mathcal{N}_d:
    \mathcal{A}_{\gamma}\rightarrow\mathcal{A}_{\gamma}$ однозначно обратим,
     его обратный   задается формулой   \eqref{sp__4_} с $\lambda=0$.

\section {Построение собственных функций оператора $\nabla div$}
\subsection {Связь между собственными функциями операторов $\nabla div$
и Лапласа-Неймана }
\begin{problem} Найти  все ненулевые собственные значения $\mu
$ и  собственные вектор-функции $\mathbf{u}(\mathbf{x})$ в
${{\mathbf{L}}_{2}}(G)$ оператора градиент дивергенции такие, что
\begin{equation}  \label{gd   1}-\nabla\mathrm{div }\mathbf{u}=
\mu \mathbf{u}\quad  \text {в} \quad G,\quad \mathbf{n}\cdot
\mathbf{u}{{|}_{\Gamma }}=0,
\end{equation}
              где
$\mathbf{n}\cdot \mathbf{u}$ - проекция вектора $\mathbf{u}$ на
нормальный вектор $\mathbf{n}$.\end{problem}
 К области определения
$\mathcal{M}_{(-\mathcal{N}_d)}$
 оператора $-\mathcal{N}_d $ задачи 2 отнесем все вектор-функции
 $\mathbf{u}(\mathbf{x})$ класса $\mathcal{C}^2(G)\cap
 \mathcal{C}(\overline{G})$, которые удовлетворяют граничному условию
 $ \mathbf{n}\cdot
\mathbf{u}{{|}_{\Gamma }}=0$ и условию $\nabla \text{ div
}\,\mathbf{u}\in
 {\mathbf{L}}_{2}(G)$.

Эта задача связана со спектральной задачей Неймана для скалярного
оператора Лапласа.

 \begin{problem} Найти все собственные значения
$\nu $ и собственные функции $g (\mathbf{x})$  оператора Лапласа
$-\Delta $ такие, что
      \begin{equation}  \label{gen   1}
              -\Delta g =\nu g \quad\text{в} \,\, G,\quad
 \mathbf{n}\cdot\nabla\,g|_{\Gamma }=0.
              \end{equation}\end{problem}
 К области определения
$\mathcal{M}_{\mathcal{N}_{(-\Delta)}}$
 оператора $\mathcal{N}_{(-\Delta)}$ задачи 3 относят все функции
 $g(\mathbf{x})$ класса $\mathcal{C}^2(G)\cap
 \mathcal{C}^1(\overline{G})$, такие что  
 $ \mathbf{n}\cdot\nabla
\,g|_{\Gamma }=0$ , $\Delta\,g\in
 {{L}}_{2}(G)$.

 Эта задача является самосопряженной \cite{vla, mi}.
Решения задач 2 и 3 принадлежат классу
$\mathcal{C}^\infty(\overline{G})$, так как $\Gamma\in
\mathcal{C}^\infty$. Имеет место
 \begin{lemma} Любому решению
$(\mu ,\mathbf{u})$ задачи 2 в области G соответствует
 решение $(\nu ,g )=(\mu, \text{div }\mathbf{u})$ задачи 3.
  Обратно, любому решению  $(\nu ,g )$ задачи 3
соответствует решение  $(\mu ,\mathbf{u})=(\nu, \nabla g)$ задачи
2.\end{lemma} Которая проверяется непосредственно.
    \subsection {Явные решения спектральной задачи Лапласа-Неймана в шаре }
   Согласно книге \cite{vla} В.С.Владимирова

    {\it собственные значения
    оператора Лапласа-Неймана $\mathcal{N}_\Delta$ в шаре  $B$   равны
    $\nu _{n,m}^{2}$,  где   $\nu
     _{n,m}^{{}}={{\alpha }_{n,m}}{{R}^{-1}}$,    $n\ge 0$,   $m\in N$, а   числа
      ${{\alpha }_{n,m}}>0$ суть нули  функций
      ${{\psi }_{n}^{\prime }}(z)$,
     производных ${{\psi }_{n}}(z)$, т.е. ${{\psi }_{n}}^{\prime }({{\alpha }_{n,m}})=0$.
Соответствующие $\nu _{n,m}^{2}$ собственные функции $g _{\kappa
}^{{}}$ имеют вид:
 \begin{equation}  \label{lan   1}g _{\kappa }^{{}}(r,\theta ,\varphi )=c{{_{\kappa }^{{}}}^{{}}}
 {{\psi }_{n}}{{(\alpha _{n,m}^{{}}r/R)}^{{}}}Y_{n}^{k}(\theta ,\varphi ),
  \end{equation}
где  $\kappa =(n, m, k)$- мультииндекс, \, $c_{\kappa
}^{{}}$-произвольные действительные постоянные,\,  $Y_{n}^{k}(\theta
,\varphi )$ - действительные сферические функции, \, $n\ge 0$,\,
$|k|\le n, \, m\in N$.}

Функции $g _{\kappa }^{{}}(x)$  принадлежат классу ${{C}^{\infty
}}(\overline{B})$  и при различных $\kappa$ ортогональны  в
${{L}_{2}}(B)$.   Система   функций $\{g _{\kappa }^{{}}\}$
полна в ${{L}_{2}}(B)$ \cite{mi}.
Нормируя их, получим  ортонормированный в ${{L}_{2}}(B)$ базис.

\subsection {Решение спектральной задачи 2 для $\nabla div$ в шаре}
Согласно лемме 2 вектор-функции ${{\mathbf{q}}_{\kappa }}(x)=
     \nabla {{g }_{\kappa }}(x)$
     являются решениями задачи 3 при ${\nu}_{n,m} ={\alpha }_{n,m}^2R^{-2}$ в
     ${{\mathbf{L}}_{2}}(B)$.
     Их компоненты $(q_r,q_\theta, q_\varphi)$ имеют вид
     \begin{equation}  \label{qom   1}\begin{array}{c}
     q _{r,\kappa }^{{}}(r,
     \theta ,\varphi )=c_{\kappa }(\alpha _{n,m}/R)
 {{\psi }_{n}^{{\prime}}}{{(\alpha _{n,m}^{{}}r/R)}^{{}}}Y_{n}^{k}
 (\theta ,\varphi ),\\
(q_{\varphi}+iq_{\theta})_{\kappa}=c_{\kappa }(1/r)
 {\psi }_{n}(\alpha _{n,m}r/R)\text{H}Y_{n}^{k}
 (\theta ,\varphi ).
  \end{array}\end{equation}
При $\kappa=(0,m,0)$ функция $Y_{0}^{0} (\theta ,\varphi )=1$,
$\text{H}Y_{0}^{0} =0$. Поэтому
\begin{equation}  \label{qomo   1}\begin{array}{l}
     q _{r,(0,m,0) }^{{}}(r)=c_{(0,m,0) }(\alpha _{0,m}/R)
 {{\psi }_{0}^{{\prime}}}{{(\alpha _{0,m}r/R)}},\\
(q_{\varphi}+iq_{\theta})_{(0,m,0)}=0.
  \end{array}\end{equation}
  Из этих формул легко выписать величины нормирующих множителей $c_{\kappa }$,
  при которых  $\left\|
{{\mathbf{q}}_{\kappa }}(x) \right\|=1$.

 Отметим, что $ {{\mathbf{q} }_{\kappa }}$ и $ {{\mathbf{q}
}_{{{\kappa }'}}}$ ортогональны при  ${\kappa }'\ne \kappa $.

Действительно, согласно формуле Гаусса-Остроградского
\begin{equation}  \label{ort   2}
       \int\limits_{B}\nabla g_{{\kappa }'}\cdot
       \nabla g_{\kappa }dx =-\int\limits_{B}
       g_{{\kappa}'}\Delta g_{\kappa }dx+\int\limits_{S}
       g_{\kappa}(n \cdot \nabla )g_{{\kappa}'}dS.
\end{equation}
 Функции ${{g }_{\kappa }}(x)$ являются  решениями
задачи  3  , они удовлетворяет  уравнению Гельмгольтца \eqref{gen 1}
 при $\nu ={\alpha}_{n,m}^2/{R}^2>0$ с краевым условием Неймана.
  Следовательно, граничный
интеграл пропадает, а
\begin{equation}  \label{ort   2}
 \int\limits_{B}{ {{\mathbf{q} }_{{{\kappa }'}}}\cdot  {{\mathbf{q}
}_{\kappa
}}^{{}}dx{{=}^{{}}}\frac{{\alpha}_{n,m}^2}{R^2}\int\limits_{B}{{{g
}_{{{\kappa }'}}}^{{}}{{g }_{\kappa }}^{{}}dx}}.\end{equation}
 Но функции
${{g}_{\kappa }}(x)$ и  ${{g }_{{{\kappa }'}}}(x)$ , согласно
(\ref{lan 1}), взаимно ортогональны  в ${{\text{L}}_{2}}(B)$ при
${\kappa }'\ne \kappa $. Значит, последний интеграл в \eqref{ort 2}
равен нулю и вектор - функции $ {{\mathbf{q} }_{\kappa }}$ и
$\mathbf{q}_{{\kappa }'}$ взаимно ортогональны в
${{\mathbf{L}}_{2}}(B)$; при этом
 $\left\|
{{\mathbf{q}}_{\kappa }}(x) \right\|=({\alpha}_{n,m}/{R})\left\|
{{{g}}_{\kappa }}(x) \right\|$.

{\bf Замечание.} Вектор-функции $\mathbf{q}_{\kappa }$ являются
 также собственными     функциями ротора с $\lambda =0$:
  $\rm{ rot}\,\mathbf{q}_{\kappa }=0$,
     $\gamma\mathbf{n}\cdot \mathbf{q}_{\kappa }=0$  .

 \subsection{Сходимость ряда Фурье по  собственным функциям
  оператора Лапласа-Неймана  в норме пространства Соболева}
В $\S\ 2.5$ главы 4 книги В.П.Михайлова \cite{mi} для областей $G$ с
границей $\Gamma\in \mathcal{C}^s$ определены подпространства
$H^s_\mathcal{N}(G)$  в $H^s(G)$:
\begin{equation} \label{sdp 1} H^s_\mathcal{N}(G)=\{f\in H^s(B):
\gamma(\mathbf{n}\cdot\nabla)\,f=0,...,\gamma(\mathbf{n}\cdot\nabla)
  \triangle^{\sigma}f=0\}, \end{equation}
где $\sigma=[s/2]-1$, а $[s/2]$ равна целой части  числа  $s/2$,
  $ s\geq 2$, и
    $H^0_\mathcal{N}(G)=L_2(G)$,
  $H^1_\mathcal{N}(G)=H^1(G)$ по определению.
Доказано, что принадлежность $f$ пространству $H^s_\mathcal{N}(G)$
необходима и достаточна для сходимости ее ряда Фурье по системе
собственных функций ${g}_{\kappa}$ оператора Лапласа-Неймана в
$H^s(G)$(см. теоремы 8 и 9 $\S\ 2.5$ гл. 4).

\subsection{Полнота системы собственных вектор-функций
 оператора градиент дивергенции в пространстве $\mathcal{A}_\gamma$}

Действительно, каждый элемент $\mathbf{q}_{\kappa }(x)=
\nabla\,g_{\kappa} $ принадлежит пространству
$\mathcal{A}_{\gamma}$, так как $g_{\kappa }\in H^1(B)$ и
$\gamma\mathbf{n}\cdot\nabla g_{\kappa }=0$. С другой стороны,
функция $h$ из $ H^1(B)$ разлагается в сходящийся в среднем ряд
\begin{equation}
\label{rh   1}h=\sum_{\kappa } (h,\widehat{g}_{\kappa
})\widehat{g}_{\kappa },  \quad \widehat{g}_{\kappa
}=({\alpha}_{n,m}/{R}){g}_{\kappa },\quad
(\widehat{g}_{\kappa},\widehat{g}_{\kappa '})= \delta_
{\kappa,{\kappa}'}.
\end{equation}
Следовательно,
\begin{equation}
\label{rnh   2}\nabla h=\sum_{\kappa } (h,\widehat{g}_{\kappa
})\nabla\widehat{g}_{\kappa }=\sum_{\kappa } (h,\widehat{g}_{\kappa
})\mathbf{q}_{\kappa }.\end{equation}

\subsection{Сходимость ряда $\mathbf{f}$ 
  по собственным функциям оператора $\nabla\mathrm{div}$ в норме пространства Соболева $H^s(B)$} Определим подпространство
   $\mathcal{A}^s_\mathcal{K}(B)$ в $\mathcal{A}$
при $s\geq 2$:
\[ \mathcal{A}^s_\mathcal{K}(B)= \{\mathbf{f}\in
\mathcal{A}\cap\mathbf{H}^s(B):
\gamma\mathbf{n}\cdot\mathbf{f}=0,...,
 \gamma\mathbf{n}\cdot\ (\nabla\text{div})
 ^{\sigma}\mathbf{f}=0\},\,\,
 \|\mathbf{f}\|_{\mathcal{A}_\mathcal{K}^s}=
\|\mathbf{f}\|_{\mathbf{H}^s},\]
 где  $\sigma=[s/2]-1$, а $[s/2]$-целая часть $s/2$;\quad
$\mathcal{A}_\mathcal{K}^0(B)=\mathbf{L}_2(B)$,\newline
  $\mathcal{A}_\mathcal{K}^1(B)=\mathbf{H}_\gamma^1(B)$ по определению.
  Имеет место
\begin{theorem}  Для того, чтобы $\mathbf{f}\in
\mathcal{A}$
 разлагалась в ряд Фурье
 \begin{equation} \label{arof 1}
\mathbf{f}(\mathbf{x})=\sum_{\kappa}
(\mathbf{f},\mathbf{q}_{\kappa})\mathbf{q}_{\kappa}(\mathbf{x})
\end{equation}
 по системе собственных вектор-функций $\mathbf{q}_{\kappa}(\mathbf{x})$
 оператора градиента дивергенции в шаре,
 сходящийся в норме
 пространства Соболева $\mathbf{H}^s(B)$, необходимо и достаточно,
  чтобы $\mathbf{f}$ принадлежала $\mathcal{A}^s_\mathcal{K}(B)$.

 Если $\mathbf{f}\in \mathcal{A}^s_\mathcal{K}(B)$,
то сходится ряд
\begin{equation} \label{arof 2}
\sum_{\kappa}{\nu}_{\kappa}^{2s}\,
|(\mathbf{f},\mathbf{q}_{\kappa})|^2 ,\quad
{\nu}_{\kappa}=({\alpha}_{n,m})/R
\end{equation} и существует такая положительная постоянная $C>0$, не
зависящая от $\mathbf{f}$, что
\begin{equation} \label{aorf 3}
\sum_{\kappa} {\nu}_{\kappa}^{2s}\,
|(\mathbf{f},\mathbf{q}_{\kappa})|^2 \leq
C\|\mathbf{f}\|^2_{\mathbf{H}^s(B)}.
\end{equation}
  Кроме того, если $s\geq 2$, то любая вектор-функция  $\mathbf{f}$ из
  $\mathcal{A}^s_\mathcal{K}(B)$
разлагается в в ряд Фурье, сходящийся в пространстве $\mathbf{C}^{s-2}(\overline{B})$.%
\end{theorem}

Действительно, граница шара $S\in \mathcal{C}^\infty$ и собственные
 функции  ${q}_{\kappa}(\mathbf{x})$  оператора градиент дивергенции
 в шаре:\,$-\nabla\text{div}\,{q}_{\kappa}(\mathbf{x})=
 {\nu}_{\kappa}^2{q}_{\kappa}(\mathbf{x})$,\,\newline
 $\gamma\mathbf{n}\cdot\mathbf{q}_{\kappa}(\mathbf{x})=0$,
 принадлежат классу
 $\mathcal{C}^{\infty}$  в  $\overline{B}$ . Значит, они
 и их конечные линейные комбинации 
 $\sum_{\kappa}
c_{\kappa}\,{q}_{\kappa}(\mathbf{x})$ принадлежат любому из
пространств  $\mathcal{A}^l_\mathcal{K}(B)$ 
 и  $\gamma\mathbf{n}\cdot(\nabla \text{div})^m\,\mathbf{q}_{\kappa}
 (\mathbf{x})=0$ при $m=0,1,...$\,. %
 и т. д.

 Согласно $\S\ 5.1 $ главы 3 в \cite{mi} для нормальной компоненты
 следа вектор-функции
 $\mathbf{f}\in  \mathbf{H}^l(B)$ на $S$ и
 ее производных
$\partial^\alpha \mathbf{f}$ при $|\alpha|<l$ имеются оценки
\begin{equation}
  \label{onaf  4}
  \|\gamma\mathbf{n}\cdot\partial^{\alpha} \mathbf{f}\|_{L_2(S)}\leq
  \|\gamma\partial^{\alpha} \mathbf{f}\|_{\mathbf{L}_2(S)}\leq c \| \mathbf{f}|\|_{\mathbf{H}^{|\alpha|+1}(B)}
  \leq c \|\mathbf{f}|\|_{\mathbf{H}^{l}(B)}.
  \end{equation}

 Обозначим через $\mathbf{S}_N(\mathbf{x})$ частичную сумму ряда
 \eqref{arof 1},   $\mathbf{S}_N(\mathbf{x})\in \mathbf{A}^s_\mathcal{K}(B)$ при всех
$N\geq 1$ и $s\geq 1$. Рассмотрим разность
$\mathbf{f}(\mathbf{x})-\mathbf{S}_N(\mathbf{x})$ и воспользуемся
оценкой \eqref{onaf 4}  следа   на $S$ ее нормальной компоненты:
\begin{equation}
  \label{oaf  0}
  \|\gamma\mathbf{n}\cdot (\mathbf{f}-\mathbf{S}_N)\|_{L_2(S)}\leq
   \|\gamma (\mathbf{f}-\mathbf{S}_N)\|_{\mathbf{L}_2(S)}\leq
  c \| (\mathbf{f}-\mathbf{S}_N)|\|_{{\mathbf{H}}^{1}(B)}
    \end{equation}
 Если
ряд Фурье \eqref{arof 1} функции $\mathbf{f}\in  \mathcal{A}$
сходится в норме $\mathbf{H}^1(B)$, то
$$\|\mathbf{f}(\mathbf{x})-\mathbf{S}_N(\mathbf{x})\|_{{\mathbf{H}}^{1}(B)}\rightarrow
0\,\, \text{при}\,\, N\rightarrow \infty.$$ Так как
$\gamma(\mathbf{n}\cdot\mathbf{S}_N)=0$ при любых $N$, то $ \|
\gamma\mathbf{n}\cdot\mathbf{f}\|_{L_2(S)}=0$ и, значит, $
\gamma\mathbf{n}\cdot\mathbf{f}=0$ и $\mathbf{f}\in
\mathbf{A}^1_\mathcal{K}(B)$.

 Если
ряд Фурье \eqref{arof 1} функции $\mathbf{f}$ сходится в норме
$\mathbf{H}^2(B)$, то он сходится и в норме $\mathbf{H}^1(B)$,
Значит, $ \gamma\mathbf{n}\cdot\mathbf{f}=0$ и $\mathbf{f}\in
\mathbf{A}^2_\mathcal{K}(B)$.

 Доказываем индукцией по $s$. Пусть $p>1$ и утверждение справедливо
 при  $s<2p$.    Покажем, что если
ряд Фурье \eqref{arof 1} функции $\mathbf{f}\in  \mathcal{A}$
сходится в норме $\mathbf{H}^s(B)$ при $s=2p$, то $\mathbf{f}\in
\mathbf{A}^{2p}_\mathcal{K}(B)$. 

Действительно, согласно оценке
\eqref{onaf 4} следа нормальной компоненты градиента дивергенции в
степени $\sigma=[s/2]-1=p-1$:
 \[ \|\gamma\mathbf{n}\cdot (\nabla\text{div})^{\sigma}\,(\mathbf{f}-\mathbf{S}_N)\|_{L_2(S)}\leq
   \|\gamma(\nabla\text{div})^{\sigma}\,
   (\mathbf{f}-\mathbf{S}_N)\|_{\mathbf{L}_2(S)}\leq
  c \| (\mathbf{f}-\mathbf{S}_N)|\|_{{\mathbf{H}}^{2p-1}(B)}.\]
Так как $\gamma(\mathbf{n}\cdot
(\nabla\text{div})^{\sigma}\,\mathbf{S}_N)=0$ при любых $N$, то
аналогично предыдущему
 $\gamma\mathbf{n}\cdot\mathbf{f}=0$,...,
 $\gamma\mathbf{n}\cdot (\nabla\text{div})^{\sigma} \mathbf{f}=0$.
 Значит, $\mathbf{f}\in
\mathbf{\Lambda}_{\mathcal{K}}^{2p}(B)$.

 Далее, если ряд Фурье
\eqref{arof 1} функции $\mathbf{f}$ сходится в норме
$\mathbf{H}^s(B)$, $s=2p+1$ то $\sigma=p-1$ и аналогично предыдущему
$\mathbf{f}\in \mathbf{A}^s_\mathcal{K}(B)$ при $s=2з+1$ .
Необходимость доказана.

Пусть $\mathbf{f}\in \mathbf{\Lambda}^s_\mathcal{K}(B)$, где $s>0$ .
Установим справедливость неравенства \eqref{aorf 3}.

 Так как $(-\nabla\text{div})\,
\mathbf{q}_{\kappa}=\nu_\kappa^2 \mathbf{q}_{\kappa}$,\,
$\gamma\mathbf{n}\cdot\mathbf{q}_\kappa=0$ \,и
\,$\gamma\mathbf{n}\cdot\mathbf{f}=0$,
 то согласно формуле Грина \eqref{gri__2_}
\begin{equation}
\label{arot   3}
(-\nabla\text{div}\,{\mathbf{f}},\mathbf{q}_{\kappa}) =
({\mathbf{f}},-\nabla\text{div}\,\mathbf{q}_{\kappa})
=\nu^2_\kappa({\mathbf{f}} ,\mathbf{q}_{\kappa}).
\end{equation}

а) Пусть $s=2p$. Обозначим через $\beta_{\kappa}$ коэффициенты Фурье
функции $(-\nabla\text{div})^p \mathbf{f}$. Согласно формуле
\eqref{arot 3}, учитывая, что $\mathbf{f}\in
\mathbf{\Lambda}^s_\mathcal{K}(B)$, имеем
\begin{equation}
\label{agd   s}\beta_{\kappa}=((-\nabla\text{div})^p\, \mathbf{f},
\mathbf{q}_{\kappa})=\nu^2_\kappa((-\nabla\text{div})^{p-1}\,{\mathbf{f}}
,\mathbf{q}_{\kappa})=...=\nu^{2p}_\kappa
({\mathbf{f}},\mathbf{q}_{\kappa}).\end{equation}

 Поскольку $(-\nabla\text{div})^p\mathbf{f}\in
\mathbf{L}_2(B)$, то
 \begin{equation}
  \label{arodf  1}\Sigma_\kappa \beta_{\kappa}^2
=\Sigma_\kappa \nu^{2s}_\kappa ({\mathbf{f}},\mathbf{q}_{\kappa})^2=
\|(-\nabla\text{div})^{s/2} \mathbf{f}\|^2.\end{equation}
неравенство \eqref{aorf 3} доказано.

б) Пусть $s=2p+1$. Функция $(-\nabla\text{div})^p \mathbf{f}\in
\mathbf{H}^1(B)$. Поэтому имеет место неравенство
 \[\Sigma_\kappa
\beta_{\kappa}^2 \nu_{\kappa}^2=\Sigma_\kappa \nu^{2(2p+1)}_\kappa
({\mathbf{f}},\mathbf{q}_{\kappa})^2\leq\, c \|(-\nabla\text{div})^p
\mathbf{f}\|^2_{\mathbf{H}^1(B)},\] из которого по Лемме 3 $\S 2.5$
гл. 4 в \cite{mi} вытекает неравенство \eqref{aorf 3}.

Вернемся к   частичной сумме $\mathbf{S}_l(\mathbf{x})$ ряда
\eqref{arof 1}. Как мы уже отмечали $\mathbf{S}_l(\mathbf{x})\in
\mathbf{A}^s_\mathcal{K}(B)$ при всех $l>0$. В частности, имеем
\quad $\text{rot} \mathbf{S}_l(\mathbf{x})=0$ и $
\gamma\mathbf{n}\cdot\nabla\text{div}\mathbf{S}_l=0$.

Поэтому оценка \eqref{apgro s} при $s=0$ принимает вид
\begin{equation} \label{apr 1} C_1\|\mathbf{S}_l\|_{2}
\leq\|\nabla\mathrm{div}\,\mathbf{S}_l\|+
\|\mathbf{S}_l\|.\end{equation}

Легко видеть, что $\|\mathbf{S}_l\|^2\leq
c\|\nabla\mathrm{div}\mathbf{S}_l\|^2$, где
$c=max_{m,n}\nu_{m,n}^{-4}$.
 Поэтому
\begin{equation} \label{apr 2} \|\mathbf{S}_l\|^2_{2 }
\leq a_1\|\nabla\mathrm{div}\,\mathbf{S}_l\|^2.\end{equation}
Следовательно, по индукции при $s=2p>2$
\begin{equation} \label{apr s} \|\mathbf{S}_l\|^2_{s}
\leq a_p\|(\nabla\mathrm{div})^p\,\mathbf{S}_l\|^2.\end{equation}
Пусть $\mathbf{f}\in \mathbf{A}^s_\mathcal{K}(B)$, где $s=2p>0$.
Согласно неравенству \eqref{aorf 3}, ряды в правой части \eqref{apr
s} сходятся и если $l>m\geq 1$, то
\[\|\mathbf{S}_l-\mathbf{S}_m\|^2_{s}\leq
 a_p(\|(\nabla\mathrm{div})^p
(\mathbf{S}_l-\mathbf{S}_m)\|^2 =\]
\[a_p\sum_{m+1}^l\nu_{\kappa}^{2p}|\mathbf{f},\mathbf{q}_{\kappa})|^2
 \rightarrow 0\]
 при
$l,m\rightarrow\infty$. Это означает, что ряд \eqref{arof 1}
сходится
к $\mathbf{f}$ в $\mathbf{H}^s(B)$.

Далее, при $s\geq 2$ в трехмерном шаре $B$ имеется вложение
пространств $\mathbf{H}^s(B)\subset\mathbf{C}^{s-2}(\overline{B})$ и
оценка:
\begin{equation}
  \label{cw  1}\|\mathbf{f}\|_{\mathbf{C}^{s-2}(\overline{B})}\leq C_s\|
\mathbf{f}\|_{s}\end{equation} для любой функции $\mathbf{f}\in
\mathbf{H}^s (B)$, в которой постоянная  $C_s>0$ не зависит от
$\mathbf{f}$ (см., например, Теорему 3 $\S\,6.2$ в \cite{mi})). В
частности,
\begin{equation}
  \label{cs  l}\|\mathbf{S}_l-\mathbf{S}_m\|_{\mathbf{C}^{s-2}(\overline{B})}
  \leq C_s\|\mathbf{S}_l-\mathbf{S}_m\|_{s}.\end{equation} Если
  $\|\mathbf{S}_l-\mathbf{S}_m\|_{s}\rightarrow 0$ при
$l,m\rightarrow\infty$, то
$\|\mathbf{S}_l-\mathbf{S}_m\|_{\mathbf{C}^{s-2}(\overline{B})}\rightarrow
0$. Это означает, что ряд \eqref{arof 1} сходится к $\mathbf{f}$ в
$\mathbf{C}^{s-2}(\overline{B})$.
Теорема доказана.

 {\bf Следствие.} {\it  Вектор-функция $f$ из
$\mathcal{A}\cap\mathbf{C}^{\infty}_0({B})$ разлагается в
 ряд Фурье \eqref{arof 1},
сходящийся в любом из пространств $\mathbf{C}^{n}(\overline{B})$,
$n\in \mathrm{N}$.}

\subsection{Скалярное произведение функций $\mathbf{f}$ и
$\mathbf{g}$ из $\mathcal{A}_{\gamma}$ в базисе из собственных
функций градиента дивергенции} Оно имеет вид:
\begin{equation} \label{spa 2}( \mathbf{f}, \mathbf{g})=
\sum_{\kappa} \,
(\mathbf{f},\mathbf{q}_{\kappa})(\mathbf{g},\mathbf{q}_{\kappa})
\end{equation}
Если $\mathbf{f}$ и $\mathbf{g}$ принадлежат
$\mathcal{A}^2_\mathcal{K}(B)$, то
\begin{equation} \label{sgd }(-\nabla\bf{div}\, \mathbf{f}, \mathbf{g})=
(\mathbf{f},-\nabla\text{div}\mathbf{g})= \sum_{\kappa}
{\nu}_{\kappa}^2[(\mathbf{f},\mathbf{q}_{\kappa})
(\mathbf{g},\mathbf{q}_{\kappa}) ].
\end{equation}Ясно, что оператор $\nabla\mathbf{div}$ является
самосопряженным в $\mathcal{A}_{\gamma}$.



\newpage
Реферат:

Оператор  градиент дивергенции  в подпространствах
$\mathbf{L}_{2}(G)$

  Р.С.Сакс

Автор изучает структуру пространства
  $\mathbf{L}_{2}(G)$ вектор-функций, квадратично интегрируемых
      по ограниченной односвязной  области $G$ трехмерного    пространства 
     с гладкой границей 
    и роль операторов   градиента дивергенции
     и ротора в построении базисов в  подпространствах
     ${\mathcal{{A}}}$ и  ${\mathcal{{B}}}$.

     Доказана       самосопряженность расширения $\mathcal{N}_d$
   оператора $\nabla\mathrm{div}$ в подпространство
   $\mathcal{A} _ {\gamma}\subset {\mathcal{{A}}}$ и базисность
   системы его   собственных функций. Выписаны явные
   формулы решения спектральной   задачи в шаре и условия разложимости
   вектор-функции   в ряд Фурье по собственным функциям
   градиента дивергенции. Изучена разрешимость краевой
   задачи:\newline
   $\nabla\mathrm{div}\,\mathbf{u}+\lambda\,\mathbf{u}=\mathbf{f}$\, в\, $G$, \,
       $ (\mathbf {n}\cdot\mathbf {u})|_{\Gamma}=g$ в пространствах
         Соболева $\mathbf{H}^{s}(G)$ порядка $s\geq 0$ и в подпространствах.

 Попутно изложены аналогичные результаты
         для оператора ротор и его симметричного расширения $S$ в $\mathcal{B}$.

Operator gradient of divergencie in subspaces of
 $\mathbf{L}_{2}(G)$ space.

  Saks Romen Semenovich

Сакс Ромэн Семенович ведущий научный сотрудник
 Институт Математики с
ВЦ УНЦ РАН 450077, г. Уфа, ул. Чернышевского, д.112 телефон: (347)
272-59-36
                 (347) 273-34-12
факс:        (347) 272-59-36 телефон дом.: (347) 273-84-69 моб.
+79173797538

 e-mail: romen-saks@yandex.ru

 \end{document}